\documentclass[12pt]{article}
\usepackage{amsfonts,amsmath,amsthm,amssymb,graphicx,mathrsfs}
\usepackage[fit]{truncate}
\usepackage{fullpage}
\usepackage{hyperref}

\theoremstyle{plain}
\newtheorem{theorem}{Theorem}
\newtheorem{lemma}[theorem]{Lemma}
\newtheorem{assume}{Assumption}
\newtheorem{remark}{Remark}

\allowdisplaybreaks

\begin{document}

\title{A Bayesian Approach for Noisy Matrix Completion:
Optimal Rate under General Sampling Distribution}

\author{The Tien Mai$ ^{(1,2)}_{\footnote{mai.thetien@insight-centre.org;
\url{http://sites.google.com/site/thetienmai/}}} $,
Pierre Alquier$ ^{(3)}_{\footnote{pierre.alquier@ensae.fr; 
\url{http://alquier.ensae.net/}}} $
\\
{\small $  ^{(1)} $ School of Mathematical Sciences, University College Dublin}
\\
{\small $  ^{(2)} $ Insight Centre for Data Analytics, Ireland}
\\
{\small $  ^{(3)} $ ENSAE-CREST}
}

\date{}

\maketitle

\begin{abstract}
Bayesian methods for low-rank matrix completion with noise have been shown to be very
efficient computationally~\cite{alquier2014bayesian, lawrence2009non, lim2007variational,
salakhutdinov2008bayesian, zhou2010nonparametric}. While the behaviour of penalized
minimization methods is well understood both from the theoretical and computational
points of view
(see~\cite{candes2010matrix, candes2010power, koltchinskii2011nuclear, recht2013parallel}
among others) in this problem, the theoretical optimality of Bayesian estimators have not
been explored yet. In this paper, we propose a Bayesian estimator for matrix completion
under general sampling distribution. We also provide an oracle inequality for this
estimator. This inequality proves that, whatever the rank of the matrix to be estimated,
our estimator reaches the minimax-optimal rate of convergence (up to a logarithmic factor).
We end the paper with a short simulation study.
\end{abstract}

\section{Introduction}
The ``Netflix Prize''~\cite{bennett2007netflix} generated a significant
interest in the \textit{matrix completion} problem. The Netflix data can be
represented as a sparse matrix made up of ratings given by users (rows) to
movies (columns). To infer the missing entries is thus very helpful to propose
sensible advertisement and improve the sales. However, it is totally impossible to
recover an uncomplete matrix without any assumption. A suitable condition, popular
in practice for this problem, is that the matrix has low-rank or approximately low-rank~\cite{
alquier2013bayesian, alquier2014bayesian, candes2010matrix, candes2009exact,
candes2010power, klopp2014noisy, koltchinskii2011nuclear}. For the Netflix problem,
this assumption is sensible as it means that many movies (or users) have similar profiles.

Let $ M_{m\times p}^0 $ be an unknown matrix (expected to be low-rank) and $ (X_1, Y_1),
\ldots, (X_n, Y_n) $ be i.i.d random variables drawn from a joint distribution
$\mathbf{P}$. We assume that
\begin{equation}
\label{main model}
Y_i = M^0_{X_i} + \mathcal{E}_i, \quad i = 1, \ldots, n,
\end{equation}
the noise variables $ \mathcal{E}_i $ are independent from $X_i$ and
$ \mathbb{E} (\mathcal{E}_i) = 0. $ We let $\Pi$ denote the marginal
distribution of $X$ when $(X,Y)\sim\mathbf{P}$. Remark that $\Pi$ is a
distribution on the set
$ \mathfrak{X} = \lbrace 1, \ldots, m \rbrace \times \lbrace 1, \ldots, p \rbrace $.
Then, the problem of estimating $ M^0 $ with $ n < mp $ is called the noisy
matrix completion problem under general sampling distribution.

A special instance of this problem is that the sampling distribution $ \Pi $  is uniform,
this assumption is done for example in~\cite{alquier2014bayesian, candes2010matrix,
candes2009exact, candes2010power, koltchinskii2011nuclear}. Clearly, in practice, the observed
entries are not always uniformly distributed: for example, some movies are more
famous than others, and thus receive much more ratings. More importantly, the sampling
distribution is not known in practice. More general sampling schemes than uniform distribution
had been already studied, see e.g.~\cite{foygel2011learning, klopp2014noisy,
negahban2012restricted}, but there are still some assumptions on $\Pi$ in these
papers. Here, we do not impose any restriction on $\Pi$. From now,
$ \Pi_{ij} = \mathbb{P} \left( X = \{i,j\} \right) $ will denote the probability to observe the
$ (i,j) $-th entry.

For any matrix $ A_{m\times p} $, let $\|A\|_F$ denote the Frobenius norm, i.e, $\|A\|_F^2
= {\rm Tr}(A^T A)$.  We define a ``generalized Frobenius norm'' as follows
\begin{equation*}
\|A\|^2_{F,\Pi} = \sum\limits_{ij} (A_{ij})^2 \Pi_{ij}.
\end{equation*} 
Note that when the sampling distribution $ \Pi $ is uniform, then $ \|A\|^2_{F,\Pi}
= (1/mp)\|A\|^2_F. $
For any matrix $ M_{m \times p} \in \mathbb{R}^{mp} $, we define the empirical risk as
\begin{equation*}
r(M) = \dfrac{1}{n} \sum\limits_{i = 1}^{n} \left( Y_i - M_{X_i} \right) ^2 
\end{equation*}
and the prediction risk 
\begin{equation*}
R(M) = \mathbb{E}_{(X,Y)\sim\mathbf{P}} \left[ \left(Y-M_X\right) ^2  \right].
\end{equation*}
In this paper, the prediction problem is considered, i.e, the objective is to define
an estimator $ \widehat{M} $ such that $ R(\widehat{M}) - R(M^0) $ is as small as possible.
Remark that $  R(M) - R(M^0) =  \| M - M^0  \|^2_{F,\Pi}$ for any $ M $
(using Pythagorean Theorem).

When handing with this problem, most of the recent methods are often based on minimizing a
criterion of the fit to the observations, such as $ r(M) $, penalized by the nuclear-norm
or the rank of the matrix. A first result can be found in by Cand{\`e}s and
Recht~\cite{candes2009exact}, Cand{\`e}s and Tao~\cite{candes2010power} for exact matrix
completion (noiseless case, i.e. $  \mathcal{E}_i = 0 $). These results were then developed
in the noisy case~\cite{candes2010matrix, koltchinskii2011nuclear}.
Some efficient algorithms had also been proposed, for example see~\cite{recht2013parallel}.

Recently, some authors have studied a more general problem, the so-called \textit{Trace
regression} problem: \cite{klopp2014noisy, koltchinskii2011nuclear}. This problem includes
matrix completion, together with other well-known problems (linear regression, reduced rank
regression and multitask learning) as special cases. They proposed nuclear-norm penalized
estimators and provided reconstruction errors for their methods. They also proved that these
errors are minimax-optimal (up to a logarithmic factor). Note that the average quadratic error
on the entries of a rank-$r$ matrix size $ m \times p $ from $ n $-observations can not be
better than: $ r\max (m,p)/n $~\cite{koltchinskii2011nuclear}.

On the other hand, Bayesian methods have been also considered~\cite{alquier2014bayesian,
lawrence2009non, lim2007variational, salakhutdinov2008bayesian, zhou2010nonparametric}.
Most Bayesian estimators are based on conjugate priors which allow to use Gibbs
sampling~\cite{alquier2014bayesian, salakhutdinov2008bayesian} or Variational Bayes
methods~\cite{lim2007variational}. These priors are discussed in details
in~\cite{alquier2014bayesian}. These algorithms are fast enough to deal
with large datasets like Netflix or
MovieLens\footnote{\url{http://grouplens.org/datasets/movielens/}}, and are actually
tested on these datasets in those papers.
However, the theoretical understanding of Bayesian algorithms is not satisfying. Up to our
knowledge, the minimax-optimality - and even the consistency - of the Bayesian
estimator under conjugate prior is an open question.

In this paper, we design a new prior and prove an minimax-optimal oracle bound
for the corresponding Bayesian estimator. This is presented in Section~\ref{section_theorem}.
In Section~\ref{section_simulation}, we discuss the implementation of our
Bayesian estimator. Some experiments comparing our estimator to the one based
on conjugate priors are done on
simulated datasets. The proof of the main result is provided in the appendix.

\section{Main Result}
\label{section_theorem}

Before we introduce our estimator, let us formulate some assumptions.

\begin{assume}
\label{bounded assume}
There is a known constant $ L $ such that
\begin{equation*}
 \| M^0 \|_{\infty} = \sup\limits_{i,j} \vert M^0_{ij} \vert   \leq  L  < +\infty.
\end{equation*}
\end{assume}
This is a mild assumption. In the Netflix and MovieLens datasets, the ratings
belong to the set $ \{ 1 , 2, 3, 4, 5 \} $, so we can take $L=5$.

\begin{assume}
\label{bruit-Pac}
The noise variables  $ \mathcal{E}_1, \ldots, \mathcal{E}_n$ are independent 
and independent of $X_1,\ldots,X_{n}$. There exist two known
constants $\sigma>0$ and $\xi>0$ such that
$$ \mathbb{E} (\mathcal{E}_{i}^{2})\leq \sigma^{2} $$
$$ \forall k\geq 3,\quad \mathbb{E} (|\mathcal{E}_{i}|^{k}) \leq \sigma^{2} k! \xi^{k-2}.$$
\end{assume}

Assumption \ref{bruit-Pac} states that the noise is sub-exponential, it includes the cases
where the noise is bounded or sub-Gaussian (and of course Gaussian), see e.g. Chapter 2
in~\cite{boucheron2013concentration}.

We now describe a prior $ \pi $ on matrices $ M_{m \times p} $ as follows.
Let $ K = \min (m,p) $ and $ \Gamma $ be a random variables taking value in the 
set $ \lbrace  \Gamma_1, \ldots,  \Gamma_K \rbrace $ with 
$ \mathbb{P} (\Gamma = \Gamma_k ) = \tau^{k-1} \left(   \frac{1-\tau}{1-\tau^K} \right)  $
where $ \Gamma_k = ( \overbrace{1, \ldots, 1}^{k \text{ times}}, \overbrace{0, \ldots, 0}^{K - k \text{ times}})  $
for some constant $ \tau \in (0,1) $ and $ k  \in \lbrace 1, \ldots, K  \rbrace $.
Now, assuming that $ \Gamma = \Gamma_k $ and a matrix $ M_{m \times p} $ is drawn as $ M = U_{m \times K}( V_{p \times K})^T $
where 
\begin{equation*}
U_{i,\ell}; V_{j,\ell}  \overset{{\rm i.i.d}}{\sim}
\begin{cases}
\mathcal{U} \left( \left[ - \delta, \delta \right] \right)  
&      \text{when } \Gamma_{k,\ell} = 1,
\\
\mathcal{U} \left( \left[ - \kappa, \kappa \right] \right)                                
&      \text{when } \Gamma_{k,\ell} = 0,
\end{cases}
  \quad \ell =  1, \ldots, K
\end{equation*}
with $ \delta = \sqrt{2L/K} $ and $ 0 \leq \kappa \leq (1/n)\sqrt{L/(10K)} $.
Note that, in this case, the entries of $ M $ satisfy:
$ \sup_{i,j} | M_{ij} | \leq  2L $. Moreover, when a matrix $M$ is drawn from
this prior, as $\kappa$ is small, most columns of $U$ and $V$ are almost
null. So the matrix $M=UV^T$ is very close to a rank-$k$ matrix.
Actually, the choice $\kappa=0$ leads to ${\rm rank}(M)\leq k$.

We are now ready to define our estimator. For any $\lambda>0$, we
consider the conditional probability measure $\hat{\rho}_\lambda$ given
by its density w.r.t. the probability measure $ \pi $:
\begin{equation}\label{formula_Gibbs}
\frac{d\hat{\rho}_{\lambda}}{d \pi }(M) = \frac{e^{-\lambda
r(M)}} {\int e^{-\lambda r} d \pi}.
\end{equation}
The aggregate $\widehat{M}_{\lambda}$ is defined as follows
\begin{equation}\label{Estimateur_Pierre}
\widehat{M}_{\lambda} = \int M \hat{\rho}_{\lambda}(d M).
\end{equation}
Note that, for $\lambda=n/(2\sigma^2)$, this corresponds exactly to the Bayesian
estimator that would be obtained for a Gaussian noise $\mathcal{E}_i \sim
\mathcal{N}(0,\sigma^2)$. However, a slightly different choice for $\lambda$,
denoted by $\lambda^*$ below,
will allow to obtain the optimality of the estimator under a wider class of noises.
For any $ x>0 $, define
\begin{equation*}
\mathcal{M}(x) = \left\lbrace M =UV^T, \text{ with } |U_{i\ell}| \leq \sqrt{\frac{x}{K}}, |V_{j\ell}| \leq \sqrt{\frac{x}{K}} \hspace*{3pt}  \right\rbrace .
\end{equation*}
and $ \mathcal{C} = [ 12 L  ( 2\xi + 3L)]  \vee  \left[ 8\sigma^{2} + 2(3L)^2 \right]. $
Hereafter, the main result is presented. We provide an oracle bound for our estimator
$ \widehat{M}_{\lambda^*} $.
\begin{theorem}
\label{MAIN} Let Assumption \ref{bounded assume} and \ref{bruit-Pac} be satisfied and take
$ \lambda = \lambda^* := \frac{n}{2\mathcal{C}} $.
Then, for any $\epsilon\in (0,1)$, with probability at least $1 - \epsilon$ and as soon as
$ n \geq \max(m,p) $, one has
\begin{align*}
\|  \widehat{M}_{\lambda^*}   - M^0\|_{F,\Pi}^2 
  \leq      \inf_{M \in \mathcal{M} (L)}   
       \Bigg\{   3   \| M - M^0\|_{F,\Pi}^2        
+    \mathscr{C}_{ L, \xi, \sigma, \tau }  \dfrac{  (m + p) {\rm rank} (M)  \log (K)   }{n}+
  \\
  +     \frac{8 \mathcal{C}  \log \left( \frac{2}{\varepsilon}  \right) }{n}           
            \Bigg\}   ,        
\end{align*}
where $   \mathscr{C}_{ L, \xi, \sigma, \tau }  $ is a (known) numerical constant
depending on $ L, \xi, \sigma $ and $ \tau $ only.
\end{theorem}

The proof of this theorem is given in the appendix. It
follows an argument called ``PAC-Bayesian
inequality". PAC-Bayesian inequalities were introduced in~\cite{STW,McA} in order to provide
empirical bounds on the prevision risk of Bayesian-type estimators. However, our proof
is closer to Catoni's works~\cite{catoni2003pac, catoni2004statistical, MR2483528}, where
it is shown how to derive powerful oracle inequalities from PAC-Bayesian bounds. This approach
has been used many times since then to prove oracle inequalities in many dimension-reduction
problems like sparse regression estimation~\cite{dalalyan2008aggregation,alquier2011pac,
alquier2013sparse}
or reduced-rank regression~\cite{alquier2013bayesian}.

The choice $ \lambda = \lambda^* $ comes from the proof of this theorem when optimizing
an upper bound on the risk $ R $, see~\eqref{choosing lambda} page~\pageref{choosing lambda}.
However, in practice, this choice may not be the best one. For example, in the experiments
done in Section 3 with Gaussian noise $ \mathcal{E}_i\sim
\mathcal{N}(0,\sigma^2) $, we take $ \lambda =
\frac{n}{4\sigma^2} $ that was shown in~\cite{dalalyan2008aggregation} to behave very
well in regression problems. Also, in practice, to take $K$ smaller than $\min(m,p)$
improves significantly the speed of the algorithm with little consequence on the
performance of the estimator~\cite{alquier2014bayesian}.

\begin{remark}
When $M^0\in\mathcal{M}(L)$, we can take $ M = M^0 $, one gets
\begin{align*}
\|  \widehat{M}_{\lambda^*}   - M^0\|_{F,\Pi}^2 
  \leq          \mathscr{C}_{ L, \xi, \sigma, \tau } 
   \dfrac{  (m + p) {\rm rank} (M^0)  \log (K)   }{n}
+        \frac{8 \mathcal{C}  \log \left( \frac{2}{\varepsilon}  \right) }{n}    .
\end{align*}

The rate $  (m + p) {\rm rank} (M^0)  \log (K)/n $ is minimax-optimal,
or at least almost minimax-optimal: a lower bound in this problem is provided by Theorems 5 and
7 in~\cite{koltchinskii2011nuclear}, it is $  (m + p) {\rm rank} (M^0)/n $. The optimality
of the $\log$ term is, to our knowledge, an open question. Note however that the upper
bound in~\cite{koltchinskii2011nuclear} is $  (m + p) {\rm rank} (M^0)
\log (m + p ) /n $. So, our bound represents a slight improvement in the case
$\min(m,p) \ll \max(m,p)$.
\end{remark}

\begin{remark}
When the sampling distribution $ \Pi $ is uniform in Theorem \ref{MAIN}, we obtain the following
oracle bound for the Frobenius norm
\begin{align*}
 \dfrac{1}{mp} \| \widehat{M}_{\lambda^*} - M^0\|_{F}^2    \hfill
  \leq        \inf_{M \in \mathcal{M} (L)}   
       \Bigg\{         \dfrac{3}{mp} \| M - M^0\|_{F}^2     
+    \mathscr{C}_{ L, \xi, \sigma, \tau }'  \dfrac{  (m + p) {\rm rank} (M)  \log (K)   }{n}+
\\
+        \frac{8 \mathcal{C}  \log \left( \frac{2}{\varepsilon}  \right) }{n}           
            \Bigg\}     .
\end{align*}
\end{remark}

Finally, we want to mention that the rate of~\cite{koltchinskii2011nuclear} is also
reached, in a work parallel to ours,
by Suzuki~\cite{Suzuki2014}, in a Bayesian framework. The main difference is that,
while~\cite{Suzuki2014} provides a rate of convergence in a more general low-rank tensor
estimation problem, his works do not bring an oracle inequality like Theorem~\ref{MAIN}
that can be used when $M^0$ is not exactly low-rank, but can be well approximated
by a low-rank matrix. Moreover, our result holds under any sampling distribution
$ \Pi $.

\section{Experiments and comparison with conjugate priors for simulated datasets }
\label{section_simulation}

\subsection{A Gibbs algorithm for $\widehat{M_{\lambda}}$}

As it has been shown in Section 2, our estimator $\widehat{M}_{\lambda^*}$
satisfies a powerful oracle inequality. However, as mentioned in the introduction,
the Bayesian estimator using conjugate priors is popular in practice
as it leads to a fast algorithm. The reason is that there is an explicit
form for the conditional posterior distribution of the $i$-th row of $U$,
$U_{i,\cdot}$, given the other rowss of $U$, $U_{-i,\cdot}$, and given $V$
(it is a multivariate normal distribution which parameters are known).
This allows to use a Gibbs sampler, with very good convergence properties.
This is described for example in~\cite{alquier2014bayesian} and the references
therein.

Here, straighforward but tedious computations lead to
$$ \hat{\rho}_{\lambda}(U_{i,\cdot}|k,U_{-i,\ell},V,\Gamma=\Gamma_k)
\propto \varphi\left[U_{i,\cdot};
\frac{2\lambda}{n} \Sigma_i \sum_{k:I_k=i} Y_k V_{J_k,\cdot}
,
\Sigma_i\right]
\prod_{\ell=1}^k \mathbf{1}_{\{|U_{i,\ell}| \leq \delta\}}
\prod_{\ell=k+1}^K \mathbf{1}_{\{|U_{i,\ell}| \leq \kappa\}}
$$
where we use the notation $X_1=(I_1,J_1)$, $\dots$, $X_n=(I_n,J_n)$,
$$(\Sigma_i)^{-1} = \frac{2\lambda}{n} \sum_{k:I_k=i} V_{J_k,\cdot}^T V_{J_k,\cdot}$$
and $\varphi(\cdot;m,V)$ is the density of the multivariate normal
distribution with mean vector $m$ and variance-covariance matrix $V$.
So, the conditional posterior
distribution of $U_{i,\cdot}$ is a truncated multivariate normal. To sample
from such a disitrubition is known as a very hard problem in general, see
for example~\cite{kotecha}. However, using the R package
{\bf tmvtnorm}~\cite{tmvtnorm},
it is possible to sample from a truncated multivariate normal fast enough
to compute our estimator on reasonnably large datasets.
Finally, instead of including the hyperparameter $k\in\{1,\dots,K\}$ in the
simulations, we simulated $K$ chains simultaneously, one
for every $k\in\{1,\dots,K\}$, and selected the realization of one of the chains
at each round using the probabilities given by~\eqref{formula_Gibbs}.

Also, note that the truncation procedure proposed by Suzuki in~\cite{Suzuki2014}
cannot be implemented, to our understanding, using this procedure, as the truncation
is done directly on the product $UV^T$ rather than on $U$ and $V$ individually.

\subsection{Experiments}

We use the notation $\widehat{M_{\lambda}}$ for our estimator, let us denote
$\hat{M}^{{\rm conjugate}}$ the estimator based on the Gaussian prior
for $U$ and $V$ with inverse Gamma variance, described in~\cite{alquier2014bayesian}
and in the aforementionned references. In order to compare both estimators,
a series of experiments were done with simulated data:
\begin{itemize}
 \item In the first series of simulations,
the data are simulated as
in~\cite{candes2010matrix,alquier2014bayesian}.
More precisely, a rank-$ 2 $ matrix
$ M^0_{m\times m} $ (so $m=p$) has been created as the product of two rank-$ 2 $ matrices,
$ M^0 = U^0_{m\times 2} (V^0_{m\times 2})^T $, where the entries of $ U^0 $ and
$ V^0 $ are i.i.d $ \mathcal{N} (0 , 20/ \sqrt{m}) $.
Only $ 20\% $ entries of the matrix $ M^0 $ are observed (using a uniform sampling).
This sampled set is then corrupted by noise as in~\eqref{main model}, where the
$  \mathcal{E}_i $ are i.i.d $ \mathcal{N} (0 , 1) $. We consider the cases
$m=100$, $m=200$, $m=500$ and $m=1000$.
\item The second series of simulations is similar to the first one, except that
the matrix $M^0$ is no longer rank $2$, but it can be well approximated by a rank
$2$ matrix:
$$M^0=U^0_{m\times 2} (V^0_{m\times 2})^T + \frac{1}{100}
(Z^0_{m\times 50})(W^0_{m\times 50})^T$$ where the entries of $Z^0$ and $W^0$
are i.i.d $ \mathcal{N} (0 , 20/ \sqrt{m}) $.
\item The third series of experiments is similar to the first one, but the
noise variables $  \mathcal{E}_i $ are now i.i.d from a uniform distribution
on $[-1,1]$. Note that, from a purely Bayesian point of view, this corresponds
to a mispecified model. However, the bound in Theorem~\ref{MAIN} is still valid
in this case.
\item Finally, the fourth series of experiments is similar to the first one,
noise variables $  \mathcal{E}_i $ are now i.i.d from a heavy-tailed distribution
(Student, with parameter $5$). This is another misspecified model, but in this
case, Theorem~\ref{MAIN} cannot be used.
\end{itemize}
The behavior of our estimator $ \widehat{M}_{\lambda} $ is computed through
the root-mean-squared error (RMSE) per entry,
$$ {\rm RMSE} = [(1/mp)
\|\widehat{M}_{\lambda}  -  M^0\|_F^2]^{1/2}   =
(1/m)\|\widehat{M}_{\lambda}  -  M^0\|_F .$$

\begin{table}[!ht]

\begin{center}

\begin{tabular}{  l | c | c  | c | c }

\hline\hline

prior   & $m = 100$          & $m = 200$     & $m = 500$ & $m = 1000$ 

\\

\hline

Uniform     &   0.535 ($\pm$0.003)   &   0.348        ($\pm$0.003) & 0.207 ($\pm$0.0001) & 0.141 ($\pm$0.0006)

\\

\hline

Gaussian                &           0.538 ($\pm$0.001)     &   0.345        ($\pm$0.001)  &  0.210 ($\pm$0.0001) & 0.146  ($\pm$0.001)

\\

\hline\hline

\end{tabular}

\end{center}

\caption{\textit{RMSEs in the first series of experiments
(low-rank matrix, Gaussian noise)}}

\label{results}

\end{table}

\begin{table}[!ht]

\begin{center}

\begin{tabular}{  l | c | c  | c | c }

\hline\hline

prior   & $m = 100$          & $m = 200$     & $m = 500$ & $m = 1000$ 

\\

\hline

Uniform     & 0.640 ($\pm$0.008)     &  0.387 ($\pm$0.001)  & 0.214 ($\pm$0.0008)  & 0.145 ($\pm$0.0002)

\\

\hline

Gaussian & 0.620 ($\pm$0.003)     &  0.385 ($\pm$0.001)  & 0.216 ($\pm$0.0003) & 0.145 ($\pm$0.001)

\\

\hline\hline

\end{tabular}

\end{center}

\caption{\textit{RMSEs in the second series of experiments
(approx. low-rank, Gaussian noise)}}

\label{results2}

\end{table}

\begin{table}[!ht]

\begin{center}

\begin{tabular}{  l | c | c  | c | c }

\hline\hline

prior   & $m = 100$          & $m = 200$     & $m = 500$ & $m = 1000$ 

\\

\hline

Uniform     &   0.328 ($\pm$0.002)   & 0.205 ($\pm$0.001)
                  & 0.120 ($\pm$0.001) & 0.084 ($\pm$0.002)

\\

\hline

Gaussian                &     0.334 ($\pm$0.003)    &   0.208 ($\pm$0.001)
                   &  0.126 ($\pm$0.003)  &   0.086 ($\pm$0.001)

\\

\hline\hline

\end{tabular}

\end{center}

\caption{\textit{RMSEs in the third series of experiments
(low-rank matrix, uniform noise)}}

\label{results3}

\end{table}

\begin{table}[!ht]

\begin{center}

\begin{tabular}{  l | c | c  | c | c }

\hline\hline

prior   & $m = 100$          & $m = 200$     & $m = 500$ & $m = 1000$ 

\\

\hline

Uniform     &  0.745 ($\pm$0.039)   &  0.567 ($\pm$0.005)
                 & 0.340 ($\pm$0.004) & 0.237 ($\pm$0.003)

\\

\hline

Gaussian    &  0.659  ($\pm$0.003)    & 0.439 ($\pm$0.001)
                & 0.268 ($\pm$0.002) & 0.186 ($\pm$0.002)

\\

\hline\hline

\end{tabular}

\end{center}

\caption{\textit{RMSEs in the fourth series of experiments
(low-rank matrix, heavy-tailed noise)}}

\label{results4}

\end{table}

The parameters are given as follows: for both $ \widehat{M}_{\lambda} $
and $\hat{M}^{{\rm conjugate}}$, the parameter $\lambda$ is set to
$n/4$, following~\cite{dalalyan2008aggregation}. Following~\cite{alquier2014bayesian}
we use for the parameters of the inverse Gamma prior in
$\hat{M}^{{\rm conjugate}}$ the values $a=1$, $b=1/100$. Finally, for
$ \widehat{M}_{\lambda} $, we used $\kappa=0$, $K=5$, $L=50$ and $\tau=1/2$ on all
the simulations apart from the heavy-tailed noise case, where we used $\tau=1/4$.
Note that a proper optimization with respect to the parameters $\tau$ and
$\lambda$ could lead to better results, for example through cross-validation.

\begin{figure}[!ht]
\centering
\includegraphics[scale=.5]{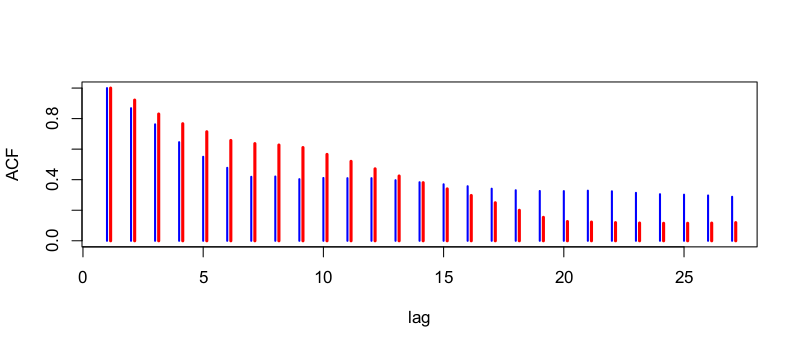}
\includegraphics[scale=.6]{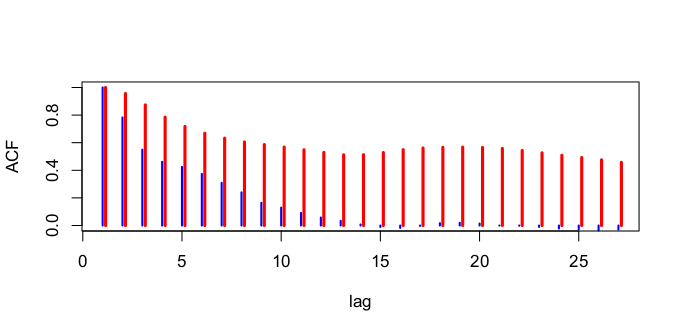}
\includegraphics[scale=.6]{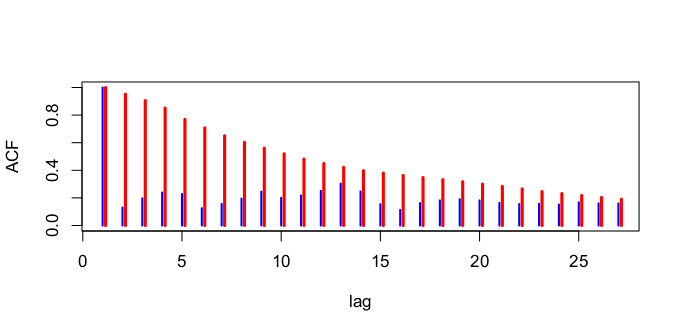}
\includegraphics[scale=.6]{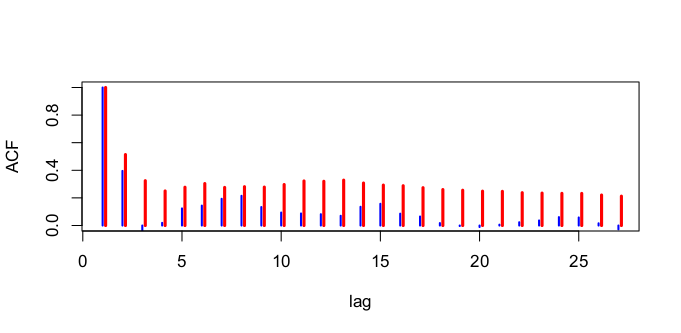}
\caption{{\small \textit{ACF of four randomly selected entries during a simulation.
These are taken from the first series of experiments. The ACF of the Gibbs sampler
for the Bayesian estimator with uniform priors, $\widehat{M_{\lambda}}$, is in
red while the ACF of the Gibbs sampler
for the Bayesian estimator with Gaussian priors, $\hat{M}^{{\rm conjugate}}$, is in
blue.}}}
\label{ACFs}
\end{figure}

The first conclusion is that the results of both methods are very close.
In many situations, however, the variance of the estimator with uniform
prior is larger than the variance of the estimator with Gaussian prior.
The evidence is that this is due to the fact that the MCMC algorithm used
to compute the estimator with Gaussian prior, $\hat{M}^{{\rm conjugate}}$,
converges faster than the algorithm used to compute the estimator with
uniform prior, $\widehat{M_{\lambda}}$. This is supported by Figure~\ref{ACFs}
page~\pageref{ACFs}. However, it seems that this difference is less and less
significant when the dimension $m$ grows.

According to our main oracle inequality, our estimator is robust to
misspecification in the low-rank assumption, see Table~\ref{results2}, and
in the noise, at least in the sub-Gaussian case, see Table~\ref{results3}.
More importantly: despite the fact that
the theoretical properties of $\hat{M}^{{\rm conjugate}}$ are not known, this
estimator is more robust than ours
to heavy-tailed noise, as shown in Table~\ref{results4}.

\section{Conclusion}

This paper proposes a Bayesian estimator for the noisy matrix completion problem
under general sampling distribution. This estimator satisfies an optimal oracle
inequality under any sampling scheme. Based on simulations, it is also clear that
this estimator performs well in practice, however, a faster algorithm for very large
datasets is still an open issue. Another important open question is the
minimax-optimality of the estimator based on Gaussian priors.

\section*{Acknowledgements}

We would like to thank the anonymous referees for their constructive comments
and Professor Taiji Suzuki for enlightening discussions.

\newpage
\section*{Appendix: Proof of Theorem \ref{MAIN}}

First, we state a version of Bernstein's inequality useful in the
proof of Theorem \ref{MAIN}. This version is taken from~\cite{MR2319879}
(Inequality 2.21 in the proof of Proposition 2.9 page 24).

\begin{lemma}
\label{lemmemassart} Let $T_{1}$, \ldots, $T_{n}$ be independent real
valued random variables. Let us assume that there are two constants
$v$ and $w$ such that
$$ \sum_{i=1}^{n} \mathbb{E}[T_{i}^{2}] \leq v $$
and for all integers $k\geq 3$,
$$ \sum_{i=1}^{n} \mathbb{E}\left[(T_{i})^{k}\right] \leq v\frac{k!w^{k-2}}{2}. $$
Then, for any $\zeta\in (0,1/w)$,
$$ \mathbb{E}
\exp\left[\zeta\sum_{i=1}^{n}\left[T_{i}-\mathbb{E}(T_{i})\right]
\right]
        \leq \exp\left(\frac{v\zeta^{2}}{2(1-w\zeta)} \right) .$$
\end{lemma}

Now, we are ready to present the proof of Theorem 1.

\begin{proof}[\textbf{Proof of Theorem \ref{MAIN}}:]

the proof is divided in two steps. In the first step, we establish a general
PAC-Bayesian inequality for matrix completion, in the style of~\cite{catoni2004statistical,
dalalyan2008aggregation}.
In the second step, we derive the oracle inequality from
the first step.

\subsubsection*{Step 1:}
Let's define, for any matrix $ M \in \mathcal{M} (2L) $, the following random variables
$$ T_{i} =   \left( Y_i - M^0_{X_{i}} \right)^{2}
                    - \left( Y_i - M_{X_{i}} \right)^{2} .$$
Note that these variables are independent. We first check that the variables $T_i$
satisfy the assumptions of Lemma~\ref{lemmemassart}, in order to apply this lemma.
We have
\begin{align*}
\sum_{i=1}^{n} \mathbb{E}[T_{i}^{2}]  &  = \sum_{i=1}^{n} \mathbb{E}
\left[
       \left( 2Y_{i} - M^0_{X_{i}} - M_{X_{i}} \right)^{2}
\left( M^0_{X_{i}} -M_{X_{i}} \right)^2
            \right]
\\
&  = \sum_{i=1}^{n} \mathbb{E} \left[
       \left( 2 \mathcal{E}_{i} + M^0_{X_i}  - M_{X_{i}}  \right)^{2}
\left( M^0_{X_{i}} -M_{X_{i}} \right)^2
            \right]
\\
&    \leq \sum_{i=1}^{n} \mathbb{E} \left[
        \left[ 8 \mathcal{E}_{i}^{2} + 2(L + 2L)^2\right]
\left[ M^0_{X_{i}} -M_{X_{i}} \right]^2
            \right]
\\
&     = \sum_{i=1}^{n} \mathbb{E} \left[ 8 \mathcal{E}_{i}^{2} + 2(3L)^2 \right]
     \mathbb{E} \left[ M^0_{X_{i}} -M_{X_{i}} \right]^2
\\
 &   \leq n \left[ 8 \sigma^{2} + 2(3L)^2\right]
\left[ R(M) - R(M^0)\right]=:v(M, M^0) = v.
\end{align*}
Next we have, for any integer $k\geq 3$, that
\begin{align*}
\sum_{i=1}^{n} \mathbb{E}\left[(T_{i})^{k}\right] \leq &
\sum_{i=1}^{n} \mathbb{E} \left[
       \left| 2Y_{i} - M^0_{X_{i}} - M_{X_{i}} \right|^{k}
\left| M^0_{X_{i}} - M_{X_{i}} \right|^k
            \right]
\\
\leq & \sum_{i=1}^{n} \mathbb{E} \left[
       2^{2k-1}\left[  |\mathcal{E}_{i}|^{k} + (L/2 + L)^{k} \right]
 \left| M^0_{X_{i}} - M_{X_{i}} \right|^{k}
            \right]
\\
\leq  &  \sum_{i=1}^{n} \mathbb{E} \left[
       2^{2k-1}\left(  |\mathcal{E}_{i}|^{k} + (\frac{3}{2} L)^{k} \right)
       (3L)^{k-2}
 \left| M^0_{X_{i}} - M_{X_{i}} \right|^{2}
            \right]
\\
\leq    &      2^{2k-1} \left[ \sigma^{2}k!\xi^{k-2}
             + \left(\frac{3}{2} L \right) ^{k} \right]  (3L)^{k-2}
\sum_{i=1}^{n}\mathbb{E}  \left| M^0_{X_{i}} - M_{X_{i}} \right|^{2}
\\
\leq  &  \frac{ \left[ \sigma^{2}k!\xi^{k-2} + ( \frac{3}{2} L )^k \right] \left[ 4 (3L) \right]^{k-2} }{ \sigma^2 + (\frac{3}{2} L)^2 } v
\\
\leq  &   \left[ k!\xi^{k-2} + \left(\frac{3}{2} L \right) ^{k-2} 
\right] [4(3 L)]^{k-2} v
\\
\leq   &   k! \left( \xi + \frac{3}{2} L \right)^{k-2}  (12L)^{k-2} v   
\leq v\frac{k!w^{k-2}}{2},
\end{align*}
with $ w:=  12 L (2 \xi + 3 L ) $.

Next, for any $\lambda\in (0,n/w)$, applying Lemma \ref{lemmemassart} with $\zeta=\lambda/n$ gives
$$
\mathbb{E} \exp\left[\lambda
\Bigl( R(M)-R(M^0)-r(M)+r(M^0)\Bigr)\right]
\leq
\exp\left[\frac{v\lambda^{2}}{2n^{2}(1-\frac{w\lambda}{n})}\right].
$$
Set $ \mathcal{C}_{\sigma, L} = 2 \left[ 4 \sigma^{2} + (3L)^2 \right] $.
For the sake of simplicity let us put
\begin{equation}
\label{defalpha}
\alpha = \left(\lambda
-\frac{\lambda^{2} \mathcal{C}_{\sigma, L}}{2n(1-\frac{w\lambda}{n})}\right) .
\end{equation}
In order to understand what follows, keep in mind that $ w $ is a constant and that our optimal estimator comes with $ \lambda = \lambda^* = \frac{n}{2\mathcal{C}} $, so $ \alpha  $ is of order $ n $.

For any $\varepsilon
>0$, the last display yields
$$
\mathbb{E} \exp\left[\alpha
                      \Bigl( R(M) - R(M^0) \Bigr)
                    +\lambda\Bigl( -r(M) + r(M^0) \Bigr)
         - \log\frac{2}{\varepsilon} \right] \leq \frac{\varepsilon}{2}.
$$
Integrating w.r.t. the probability distribution $ \pi(.) $, we get
$$
\int \mathbb{E} \exp\Biggl[\alpha
                      \Bigl( R(M) - R(M^0) \Bigr)
                    +\lambda\Bigl( -r(M) + r(M^0) \Bigr)
         - \log\frac{2}{\varepsilon}\Biggr]  \pi (d M) \leq \frac{\varepsilon}{2}.
$$
Next, Fubini's theorem gives
\begin{align*}
\mathbb{E} \int \exp\Biggl[\alpha
                      \Bigl( R(M) - R(M^0) \Bigr)
                    +\lambda\Bigl( -r(M) + r(M^0) \Bigr)
         - \log\frac{2}{\varepsilon}\Biggr]  \pi (d M)   \hspace*{1.5cm}
\\
 =  \mathbb{E} \int \exp  \left\lbrace   \alpha    \Bigl( R(M) - R(M^0) \Bigr)
                    +\lambda\Bigl( -r(M) + r(M^0) \Bigr)      \right.  -
 \\
 \left.             - \log \left[\frac{d\hat{\rho}_{\lambda}}{d \pi} (M)  \right]
        - \log\frac{2}{\varepsilon}
        \right\rbrace
         \hat{\rho}_{\lambda}(d M)
\leq \frac{\varepsilon}{2}.
\end{align*}
Jensen's inequality yields
$$
\mathbb{E} \exp\Biggl[\alpha
                      \left( \int R d \hat{\rho}_{\lambda} - R(M^0) \right)
                    +\lambda\left( -\int r d \hat{\rho}_{\lambda} + r(M^0) \right)
                 - \mathcal{K}(\hat{\rho}_{\lambda}, \pi)
         - \log\frac{2}{\varepsilon}\Biggr] \leq \frac{\varepsilon}{2},
$$
where $  \mathcal{K}(p, q) $ is the Kullback--Leibler divergence of $ p $ from $ q $.
Now, using the basic inequality $\exp(x) \geq
\mathbf{1}_{\mathbb{R}_{+}}(x)$, we get
$$
\mathbb{P}\Biggl\{ \Biggr[ \alpha
                      \left(\int R d\hat{\rho}_{\lambda} - R(M^0)\right)
                    +\lambda\left(-\int r d\hat{\rho}_{\lambda} + r(M^0) \right)
                 - \mathcal{K}(\hat{\rho}_{\lambda}, \pi)
         - \log\frac{2}{\varepsilon}\Biggr] \geq 0
\Biggr\} \leq \frac{\varepsilon}{2}.
$$
Using Jensen's inequality again gives
$$ \int R d\hat{\rho}_{\lambda} \geq R\left( \int M \hat{\rho}_{\lambda}(d M)\right)
                                      = R(\widehat{M}_{\lambda}).$$
Combining the last two displays we obtain
\begin{equation*}
\mathbb{P}\Biggl\{ R(\widehat{M}_{\lambda})  - R(M^0)
\leq 
\frac{ \int r d\hat{\rho}_{\lambda} - r(M^0) +
\frac{1}{\lambda}\left[\mathcal{K}(\hat{\rho}_{\lambda}, \pi) 
+ \log\frac{2}{\varepsilon}\right] } {\frac{\alpha}{\lambda} }
\Biggr\}
 \geq 
 1-\frac{\varepsilon}{2}.
\end{equation*}
\\
Using Donsker and Varadhan's variational inequality (Lemma 1.1.3 in Catoni~\cite{MR2483528}),
we get
\begin{equation}\label{interm3bis}
\mathbb{P}\Biggl\{ R(\widehat{M}_{\lambda}) - R(M^0)
\leq \inf_{\rho \in \mathfrak{M}_{+}^{1}(M) } \frac{ \int r
d\rho - r(M^0) +
\frac{1}{\lambda}\left[\mathcal{K}(\rho, \pi)
+ \log\frac{2}{\varepsilon}\right] } {\frac{\alpha}{\lambda} } \Biggr\} \geq 1-\frac{\varepsilon}{2},
\end{equation}
where $\mathfrak{M}_{+}^{1}(M)$ is the set of all positive
probability measures over the set of $m\times p$ matrices equiped with the Borel
$\sigma$-algebra.
\\

We now want to bound from above $ r(M) - r(M^0)$ by
$R(M)- R(M^0)$. We can use Lemma~\ref{lemmemassart} again, to
$\tilde{T}_i(\theta) = - T_i(\theta)$ and similar computations yield successively
$$
\mathbb{E} \exp\left[\lambda \Bigl( R(M^0) - R(M) +
r(M) - r(M^0) \Bigr)\right] \leq
\exp\left[\frac{v\lambda^{2}}{2n^{2}(1-\frac{w\lambda}{n})}\right],
$$
and so for any (data-dependent) $\rho$,
$$
\mathbb{E} \exp\Biggl[\beta
               \left(-\int Rd\rho + R(M^0) \right)
+ \lambda \left( \int r d\rho - r(M^0) \right) -
\mathcal{K}(\rho, \pi) - \log \frac{2}{\varepsilon}\Biggr] \leq
\frac{\varepsilon}{2},
$$
where
\begin{equation}
\label{defbeta}
\beta = \left(\lambda
+\frac{\lambda^{2} \mathcal{C}_{\sigma, L} }{2n(1-\frac{w\lambda}{n})}\right) .
\end{equation}
Here again, with the same spirit with $ \alpha $ in (\ref{defalpha}), $ \beta $ is of order $ n $ also.
So:
\begin{equation}
\label{interm4} \mathbb{P}\Biggl\{ \int rd\rho - r(M^0)
\leq \frac{\beta}{\lambda} \left[\int
Rd\rho - R(M^0) \right] + \frac{1}{\lambda}\left[
\mathcal{K}(\rho, \pi) + \log \frac{2}{\varepsilon} \right]
\Biggr\}\geq 1 - \frac{\varepsilon}{2}.
\end{equation}
\\
Combining (\ref{interm4}) and (\ref{interm3bis}) with a union bound
argument gives the general PAC-Bayesian bound
\begin{equation}
\label{PAC-bound}
\mathbb{P}\Biggl\{ R(\widehat{M}_{\lambda}) - R(M^0)
\leq 
\inf_{\rho \in \mathfrak{M}_{+}^{1}(M) } \frac{ \beta \left[\int Rd\rho -
R(M^0) \right] + 2 \left[
\mathcal{K}(\rho, \pi) + \log \frac{2}{\varepsilon} \right] } {
\alpha  } \Biggr\} 
\geq 
1-\varepsilon.
\end{equation}

\subsubsection*{Step 2:}
In the second step, we derive an explicit form for the upper bound in~\eqref{PAC-bound}.
The idea is that, if we restrict the infimum in the upper bound in~\eqref{PAC-bound}
to a small set of measures $\rho$, we are able to provide an explicit bound for this
infimum. This trick was introduced in~\cite{catoni2004statistical}.

Let $ M \in \mathcal{M} (L) $, it means that $ M = UV^T $ with 
$ |U_{i\ell}| \leq \sqrt{L/K} , |V_{j\ell}| \leq \sqrt{L/K} $.
Let us take, for any $ c $ such that $ \kappa \leq  c <  (\sqrt{2} - 1) \sqrt{L/K}$,
 the probability distribution 
$$ \rho_{U,V,c}({\rm d}\mu,{\rm d}\nu)    \propto   
\mathbf{1} (\| \mu - U \|_{\infty} \leq c, \|\nu - V \|_{\infty} \leq c) 
                             \hspace*{3pt}   \pi({\rm d}\mu,{\rm d}\nu) .$$
Note that, as $  c <  (\sqrt{2} - 1) \sqrt{L/K} $, we have 
$ {\rm supp}(\rho_{U,V,c}) \subset {\rm supp}(\pi) $ and so 
$  \mathcal{K} (\rho_{U,V,c}, \pi)  < \infty .$
\\
Thus, (\ref{PAC-bound}) becomes
\begin{align}
\label{stepDT}
\nonumber
\hspace*{-5pt}  \mathbb{P}\Biggl\{ R(\widehat{M}_{\lambda}) - R(M^0)
\leq \inf_{U,V,c} 
                   \frac{ \beta \left[\int R  d\rho_{U,V,c}   -   R(M^0) \right] 
     + 2 \left[  \mathcal{K} (\rho_{U,V,c}, \pi) + \log \frac{2}{\varepsilon} \right] } 
     {\alpha  } \Biggr\} 
\\
\geq 1-\varepsilon.
\end{align}

Let us fix $ c, U, V $. The end the proof consists in calculations to derive
an upper bound for the two terms in~\eqref{stepDT}. Firstly
\begin{align*}
\int R(M) d  \rho_{U,V,c} -    R(M^0)
&   =    \int \| \mu\nu^T - M^0 \|_{F, \Pi}^2 \hspace*{5pt} \rho_{U,V,c}({\rm d}\mu,{\rm d}\nu) 
 \\
 &  =  \int \| \mu\nu^T - U\nu^T + U\nu^T - UV^T
 + UV^T - M^0 \|_{F, \Pi}^2 \hspace*{5pt}   \rho_{U,V,c}({\rm d}\mu,{\rm d}\nu)
\\
 &  = \int \biggl( \| \mu\nu^T - U\nu^T \|_{F, \Pi}^2 + \| U\nu^T - UV^T\|_{F, \Pi}^2  +
 \\
  & + \| UV^T - M^0\|_{F, \Pi}^2
  + 2\left< \mu\nu^T - U\nu^T, U\nu^T - UV^T \right>_{F, \Pi}
  \\
  &   + 2\left< \mu\nu^T - U\nu^T, UV^T - M^0 \right>_{F, \Pi}
    \\
  &   + 2\left< U\nu^T - UV^T, UV^T - M^0 \right>_{F, \Pi}
  \biggr) \rho_{U,V,c}({\rm d}\mu,{\rm d}\nu)   .
\end{align*}
(note that we use the notation $ \left< A,B \right>_{F, \Pi} = \sum_{i,j} A_{ij} B_{ij}
\Pi_{ij} $).
As $ \int \mu \rho_{U,V,c}({\rm d}\mu) = U $ and $\int \nu \rho_{U,V,c}({\rm d}\nu) = V $,
it can be seen that integral of the three scalar products in the previous
equation vanish. Moreover,
\begin{align*}
\| (\mu - U)\nu^T \|_{F, \Pi}^2   &  = \sum\limits_{ij} \left[ (\mu - U)\nu^T \right]^2_{ij} \Pi_{ij}   
 \leq    \left(   \sup\limits_{ij}  \left[ (\mu - U)\nu^T \right]_{ij}   \right)^2 \sum\limits_{ij} \Pi_{ij}
\\
& \leq      \left(   \sup\limits_{ij}  \sum\limits_{\ell = 1}^K |\mu - U|_{i\ell}    |\nu|_{j\ell}    \right)^2
 \leq      \left(  K  \sup\limits_{i\ell} |\mu - U|_{i\ell} \,\,  \sup\limits_{j\ell}  |\nu|_{j\ell}    \right)^2
 \\
 &  \leq      \left[  K c \left(  c + \sqrt{ \dfrac{L}{K}} \right)  \right]^2  = Kc^2 ( \sqrt{K} c + \sqrt{L})^2   ,
\end{align*}
similarly $ \| U\nu^T - UV^T\|_{F, \Pi}^2     \leq       KL c ^2 $.
Therefore, from~\eqref{stepDT}, we have
\begin{align}
\label{step1stterm}
 \int \| \mu\nu^T   - M^0 \|_{F, \Pi}^2   \hspace*{4pt}  \rho_{U,V,c}({\rm d}\mu,{\rm d}\nu)      
 \leq  K c^2  \left[  ( \sqrt{K} c + \sqrt{L} )^2 +L \right]    +   \| UV^T - M^0\|_{F, \Pi}^2 .
\end{align}

So, we have an upper bound for the first term in~\eqref{stepDT}. We now deal with
the Kullback-Leibler term:
\begin{align}
\nonumber
\mathcal{K}(\rho_{U,V,c},\pi)
 =  &   \log \frac{1}{\pi(\{\mu,\nu: \|\mu - U\|_{\infty} \leq c, \|\nu - V\|_{\infty} \leq c\})}
\\
\nonumber
=   &   \log \frac{1}{\pi(\{\mu : \|\mu - U\|_{\infty} \leq c  \})}    +   \log \frac{1}{\pi(\{\nu:  \|\nu - V\|_{\infty} \leq c\})}
 \\
 \nonumber
=    &   \log \frac{1}{\int \pi(\{ \|\mu - U\|_{\infty}\leq c \}|\Gamma)\pi(\Gamma){\rm d}\Gamma} +
\\
&     \hspace*{3.6cm}  + \log \frac{1}{\int \pi(\{ \|\nu - V\|_{\infty} \leq c    \label{stepKL}
 \}|\Gamma)\pi(\Gamma){\rm d}\Gamma}.
\end{align}
\\
Note that, up to a reordering of the columns of $ U $ and $ V $, we can assume that
$U=(U_1 | \dots | U_{k_0} | 0 | \ldots | 0 | )$ and $V=(V_1 | \dots |V_{k_0}  | 0 | \ldots | 0 | )$,
where $ k_0 = {\rm rank}(UV^T)  \leq K $.
Then
\begin{align}
\nonumber
\int \pi(\{ \|\mu - U\|_{\infty}\leq c \}|\Gamma) \pi(\Gamma)     {\rm d}\Gamma     
    = \tau^{k_0 - 1}   \left(   \frac{1-\tau}{1-\tau^K} \right)    \pi(\{ \|\mu - U\|_{\infty}\leq  c \}  |\Gamma = \Gamma_{k_0})
\end{align}
and, as $\kappa\leq c$,
\begin{align*}
\pi(\{\|\mu  -  U \|_{\infty}  \leq  c \}  |\Gamma = \Gamma_{k_0})
  &    \geq 
    \prod\limits_{i=1}^m     \prod\limits_{\ell=1}^{k_0}  
     \pi(\{ |\mu_{i\ell} - U_{i\ell} |     \leq c \}  |\Gamma = \Gamma_{k_0})  \hspace*{-8pt}
      \prod\limits_{\ell=k_0 + 1}^{K}   \hspace*{-10pt}
              \pi(\{ |\mu_{i\ell} |   \leq c \}  |\Gamma = \Gamma_{k_0})
\\
&      \geq    \left( c \sqrt{\dfrac{K}{2L}} \right)^{mk_0}.
\end{align*}
\\
So,
\begin{align} \label{stepKL1}
\nonumber
\log \frac{1}{\int \pi(\{ \|\mu - U\|_{\infty} \leq c \}|\Gamma)\pi(\Gamma){\rm d}\Gamma}  
   \leq    
    ( k_0 -1) \log(1/ \tau)  + \log \left(   \frac{1-\tau^K}{1-\tau} \right)  +
  mk_0 \log \left( \dfrac{1}{c} \sqrt{\dfrac{2L}{K}} \right)
\\
     \leq    ( k_0 -1 )  \log(1/ \tau)  + \log \left(   \frac{1}{1-\tau} \right)  +
 mk_0 \log \left( \dfrac{1}{c} \sqrt{\dfrac{2L}{K}} \right)   .
\end{align}
\\
By symmetry,
\begin{align} \label{stepKL2}
  \nonumber
\log \frac{1}{\int \pi(\{ \|\nu - V \|_{\infty}\leq c \}|\Gamma)\pi(\Gamma){\rm d}\Gamma} 
  \leq        ( k_0 -1 )   \log(1/ \tau)  + \log \left(   \frac{1}{1-\tau} \right)  + \hspace*{10pt}
  \\
+  pk_0 \log \left( \dfrac{1}{c} \sqrt{\dfrac{2L}{K}} \right)   .
\end{align}
Plugging~\eqref{stepKL1} and~\eqref{stepKL2} into~\eqref{stepKL},
we obtain finally our upper bound for the Kullback-Leibler term:
\begin{align}
\label{stepKLend}
\nonumber
\mathcal{K}(\rho_{U,V,c},\pi)
&  \leq    2  ( k_0 -1 )   \log(1/ \tau)  +   2\log \left(   \frac{1}{1-\tau} \right)  +
  (m+p) k_0 \log \left( \dfrac{1}{c} \sqrt{\dfrac{2L}{K}} \right)  
  \\
&   \leq    2  k_0    \log(1/ \tau)  +   2\log \left(   \frac{\tau}{1-\tau} \right)  +
  (m+p) k_0 \log \left( \dfrac{1}{c} \sqrt{\dfrac{2L}{K}} \right)     .
\end{align}
Finally, substituting~\eqref{step1stterm} and~\eqref{stepKLend} into~\eqref{stepDT},
\begin{align*}
\nonumber
 \mathbb{P}\Biggl\{ R(\widehat{M}) -  R(M^0)
  \leq          \inf_{\begin{array}{c}U,V,c \\
U_j,V_j = 0 \text{ when } j > k_0 \end{array}} 
                   \frac{ 1 }{\alpha  } 
\Bigg[ \beta  \left(  K c^2  \left[  ( \sqrt{K} c+ \sqrt{L} )^2 +L \right] +   \right. \hspace*{1cm}
\\
+   \left.    \| UV^T - M^0\|_{F,\Pi}^2 \right)   
  +  2(m + p) k_0 \log \left( \dfrac{1}{c} \sqrt{\dfrac{2L}{K}}   \right)   +
\\
+ 4  k_0 \log(1/ \tau)  +   4\log \left(   \frac{\tau}{1-\tau} \right)     +   2\log \frac{2}{\varepsilon} 
            \Bigg]          \Biggr\}      
                   \geq 1-\varepsilon.
\end{align*}
Let us put $ c = \sqrt{(m+p)L/(18nK)}$. Note that as 
$ n \geq \max(m,p) $ then $\sqrt{(m+p)/(3n)} <1 $ and 
thus the condition $ c <   (\sqrt{2} - 1) \sqrt{L/K} $ is satisfied.
So we have the
following inequality with probability at least $ 1- \varepsilon $:
\begin{align*}
R(\widehat{M}_{\lambda}) - R(M^0) \hfill
  \leq   \hspace*{-10pt}   \inf_{\begin{array}{c} U,V \\
U_j,V_j = 0 \text{ when } j > k_0 \end{array}} \hspace*{-15pt} 
             \frac{1}{1-\frac{\lambda \mathcal{C}_{\sigma, L}}{2(n-w\lambda)}} 
\Bigg\{         \left( 1+\frac{\lambda  \mathcal{C}_{\sigma, L}}{2(n-w\lambda)}\right)
       \Bigg[      \| UV^T - M^0\|_{F.\Pi}^2 +   
\\
      +      L \dfrac{m+p}{18n} \left(  2L \dfrac{m+p}{18n} + 3L  \right)    \Bigg]
  +  \dfrac{2}{\lambda}     
\Bigg[          (m + p) k_0  \log \left( \sqrt{\frac{36n}{m+p}}  \right) +
\\
   + 2  k_0  \log(1/ \tau)  +   2\log \left(   \frac{\tau}{1-\tau} \right)     +   \log \frac{2}{\varepsilon} 
          \Bigg]
            \Bigg\}    ,
\end{align*}
where $ \alpha $ and $ \beta $ have been replaced by their definitions,
see~\eqref{defalpha} and~\eqref{defbeta}.
Taking now $ \lambda =\lambda^*= n/(2\mathcal{C}) $ with
$\mathcal{C} = \mathcal{C}_{\sigma, L} \vee   w$ in the last above display, gives 
\begin{align} \label{choosing lambda}
\nonumber
 \mathbb{P}\Biggl\{ R(\widehat{M}_{\lambda^*}) - R(M^0) \hfill
  \leq     \inf_{M \in \mathcal{M} (L)}   
       \Bigg\{       3
       \Bigg[         L^2 \dfrac{m+p }{18n} \left(  \dfrac{m+p}{9n} + 3  \right) 
        +   \| M - M^0\|_{F,\Pi}^2      \Bigg] +
 \\  
\nonumber  
 +  \dfrac{8 \mathcal{C}}{n}     
\Bigg[         \frac{1}{2} (m + p) {\rm rank (M)}  \log \left( \frac{36n}{m+p }  \right)     
                                               +  \log \frac{2}{\varepsilon}   +
\\
  + 2 {\rm rank}(M)  \log(1/ \tau)  +   2\log \left(   \frac{\tau}{1-\tau} \right)     
          \Bigg]
            \Bigg\}   
                  \Biggr\}          \geq 1-\varepsilon  ,
\end{align}
where we have used that $1-\frac{\lambda  \mathcal{C}_{\sigma, L} }{2(n-w\lambda)} \geq 1/2$
and $ 1+\frac{\lambda \mathcal{C}_{\sigma, L}}{2(n-w\lambda)}\leq 3/2$. As
$$
 \log \left( \frac{36n}{m+p}  \right)  \leq   \log \left( \frac{36mp}{ \max(m,p) }  \right)
   =   \log \left( \frac{36\min(m,p) \max(m,p)}{ \max(m,p) }  \right)  
 =  \log \left( 36K \right),
$$
we have
\begin{align}\label{last inequality}
\nonumber
 \mathbb{P}\Biggl\{ R(\widehat{M}_{\lambda^*}) - R(M^0) \hfill
  \leq     \inf_{M \in \mathcal{M} (L)}   
       \Bigg\{       3
       \Bigg[         L^2 \dfrac{m+p }{18n} \left(  \dfrac{m+p}{9n} + 3  \right) 
        +   \| M - M^0\|_{F,\Pi}^2      \Bigg] +
 \\  
\nonumber  
 +  \dfrac{8 \mathcal{C}}{n}     
\Bigg[         \frac{1}{2} (m + p) {\rm rank (M)}  \log (36K)     
                                               +  \log \frac{2}{\varepsilon}   +
\\
  + 2 {\rm rank}(M)  \log(1/ \tau)  +   2\log \left(   \frac{1}{1-\tau} \right)     
          \Bigg]
            \Bigg\}   
                  \Biggr\}          \geq 1-\varepsilon  .
\end{align}
\\
Moreover,
\begin{align*}
       L^2 \dfrac{m+p }{6n} \left(   \dfrac{m+p }{9n} + 3  \right)  
  \leq  
\mathscr{C}(L)  \dfrac{(m+p){\rm rank (M) \log (K) } }{n} ,
\end{align*}
for some constant $ \mathscr{C}(L)  >0 $ depending on $ L $ only.
Remind that $ \tau $ is a constant in $ (0,1) $, we have
\begin{align*}
2 {\rm rank}(M)  \log(1/ \tau)  +   2\log \left(   \frac{\tau}{1-\tau} \right)   
 \leq 
\mathscr{C}(\tau)    \dfrac{(m+p){\rm rank (M) \log (K) } }{n}  ,
\end{align*}
for some constant $ \mathscr{C}(\tau)   > 0 $ depending on $ \tau $ only.
Finally, from~\eqref{last inequality}, we obtain
\begin{align*}
 \mathbb{P}\Biggl\{ R(\widehat{M}_{\lambda^*}) - R(M^0) \hfill
  \leq      \inf_{M \in \mathcal{M} (L)}   
       \Bigg[   3   \| M - M^0\|_{F,\Pi}^2      +  
       \mathscr{C}(L,\mathcal{C}, \tau)  \dfrac{(m+p){\rm rank (M) \log (K) } }{n} +
\\
      +    \frac{8 \mathcal{C}  \log \left( \frac{2}{\varepsilon}  \right) }{n}    
            \Bigg] 
                  \Biggr\}          \geq 1-\varepsilon  ,
\end{align*}
for some constant $  \mathscr{C}(L,\mathcal{C}, \tau)  >0 $ depending only on
$ L, \tau $ and $ \mathcal{C} $. 
However, as the constant $ \mathcal{C} $ also depends on $ L, \xi, \sigma $ then 
$  \mathscr{C}(L,\mathcal{C}, \tau) $ can be rewritten as $  \mathscr{C}_{ L, \xi,
\sigma, \tau }$ as in the statement of the theorem.
\end{proof}

\newpage
\bibliographystyle{abbrv}

\begin{thebibliography}{10}

\bibitem{alquier2013bayesian}
P.~Alquier.
\newblock Bayesian methods for low-rank matrix estimation: short survey and
  theoretical study.
\newblock In {\em Algorithmic Learning Theory 2013}, pages 309--323. Springer,
  2013.

\bibitem{alquier2013sparse}
P.~Alquier and G.~Biau.
\newblock Sparse single-index model.
\newblock {\em The Journal of Machine Learning Research}, 14(1):243--280, 2013.

\bibitem{alquier2014bayesian}
P.~Alquier, V.~Cottet, N.~Chopin, and J.~Rousseau.
\newblock Bayesian matrix completion: prior specification.
\newblock {\em arXiv preprint arXiv:1406.1440}, 2014.

\bibitem{alquier2011pac}
P.~Alquier and K.~Lounici.
\newblock Pac-{B}ayesian bounds for sparse regression estimation with
  exponential weights.
\newblock {\em Electronic Journal of Statistics}, 5:127--145, 2011.

\bibitem{bennett2007netflix}
J.~Bennett and S.~Lanning.
\newblock The netflix prize.
\newblock In {\em Proceedings of KDD cup and workshop}, volume 2007, page~35,
  2007.

\bibitem{boucheron2013concentration}
S.~Boucheron, G.~Lugosi, and P.~Massart.
\newblock {\em Concentration inequalities: A nonasymptotic theory of
  independence}.
\newblock Oxford University Press, 2013.

\bibitem{candes2010matrix}
E.~J. Cand{\`e}s and Y.~Plan.
\newblock Matrix completion with noise.
\newblock {\em Proceedings of the IEEE}, 98(6):925--936, 2010.

\bibitem{candes2009exact}
E.~J. Cand{\`e}s and B.~Recht.
\newblock Exact matrix completion via convex optimization.
\newblock {\em Found. Comput. Math.}, 9(6):717--772, 2009.

\bibitem{candes2010power}
E.~J. Cand{\`e}s and T.~Tao.
\newblock The power of convex relaxation: near-optimal matrix completion.
\newblock {\em IEEE Trans. Inform. Theory}, 56(5):2053--2080, 2010.

\bibitem{catoni2003pac}
O.~Catoni.
\newblock {\em A PAC-Bayesian approach to adaptive classification}.
\newblock Preprint Laboratoire de Probabilit\'es et Mod\`eles Al\'eatoires
  PMA-840, 2003.

\bibitem{catoni2004statistical}
O.~Catoni.
\newblock {\em Statistical Learning Theory and Stochastic Optimization}.
\newblock Saint-Flour Summer School on Probability Theory 2001 (Jean Picard
  ed.), Lecture Notes in Mathematics. Springer, 2004.

\bibitem{MR2483528}
O.~Catoni.
\newblock {\em PAC-{B}ayesian supervised classification: the thermodynamics of
  statistical learning}.
\newblock Institute of Mathematical Statistics Lecture Notes---Monograph
  Series, 56. Institute of Mathematical Statistics, Beachwood, OH, 2007.

\bibitem{dalalyan2008aggregation}
A.~Dalalyan and A.~B. Tsybakov.
\newblock Aggregation by exponential weighting, sharp pac-bayesian bounds and
  sparsity.
\newblock {\em Machine Learning}, 72(1-2):39--61, 2008.

\bibitem{foygel2011learning}
R.~Foygel, O.~Shamir, N.~Srebro, and R.~Salakhutdinov.
\newblock Learning with the weighted trace-norm under arbitrary sampling
  distributions.
\newblock In {\em Advances in Neural Information Processing Systems}, pages
  2133--2141, 2011.

\bibitem{klopp2014noisy}
O.~Klopp.
\newblock Noisy low-rank matrix completion with general sampling distribution.
\newblock {\em Bernoulli}, 20(1):282--303, 2014.

\bibitem{koltchinskii2011nuclear}
V.~Koltchinskii, K.~Lounici, and A.~B. Tsybakov.
\newblock Nuclear-norm penalization and optimal rates for noisy low-rank matrix
  completion.
\newblock {\em The Annals of Statistics}, 39(5):2302--2329, 2011.

\bibitem{kotecha}
J. H. Kotecha and P. M. Djuric.
\newblock Gibbs Sampling Approach For Generation of Truncated Multivariate Gaussian
Random Variables.
\newblock {\em Proceedings of the IEEE Conference on Acoustics, Speech, and
Signal Processing}, 3:1757--1760, 1999.

\bibitem{lawrence2009non}
N.~D. Lawrence and R.~Urtasun.
\newblock Non-linear matrix factorization with gaussian processes.
\newblock In {\em Proceedings of the 26th Annual International Conference on
  Machine Learning}, pages 601--608. ACM, 2009.

\bibitem{lim2007variational}
Y.~J. Lim and Y.~W. Teh.
\newblock Variational bayesian approach to movie rating prediction.
\newblock In {\em Proceedings of KDD Cup and Workshop}, volume~7, pages 15--21,
  2007.

\bibitem{MR2319879}
P.~Massart.
\newblock {\em Concentration inequalities and model selection}, volume 1896 of
  {\em Lecture Notes in Mathematics}.
\newblock Springer, Berlin, 2007.
\newblock Lectures from the 33rd Summer School on Probability Theory held in
  Saint-Flour, July 6--23, 2003, Edited by Jean Picard.

\bibitem{McA}
D.~McAllester.
\newblock Some {PAC}-{B}ayesian theorems.
\newblock In {\em Proceedings of the Eleventh Annual Conference on
  Computational Learning Theory}, pages 230--234, New York, 1998. ACM.

\bibitem{negahban2012restricted}
S.~Negahban and M.~J. Wainwright.
\newblock Restricted strong convexity and weighted matrix completion: Optimal
  bounds with noise.
\newblock {\em The Journal of Machine Learning Research}, 13(1):1665--1697,
  2012.

\bibitem{recht2013parallel}
B.~Recht and C.~R{\'e}.
\newblock Parallel stochastic gradient algorithms for large-scale matrix
  completion.
\newblock {\em Mathematical Programming Computation}, 5(2):201--226, 2013.

\bibitem{salakhutdinov2008bayesian}
R.~Salakhutdinov and A.~Mnih.
\newblock Bayesian probabilistic matrix factorization using markov chain monte
  carlo.
\newblock In {\em Proceedings of the 25th international conference on Machine
  learning}, pages 880--887. ACM, 2008.

\bibitem{STW}
J.~Shawe-Taylor and R.~Williamson.
\newblock A {PAC} analysis of a {B}ayes estimator.
\newblock In {\em Proceedings of the Tenth Annual Conference on Computational
  Learning Theory}, pages 2--9, New York, 1997. ACM.

\bibitem{Suzuki2014}
T.~Suzuki.
\newblock Convergence rate of bayesian tensor estimation: optimal rate without
  restricted strong convexity.
\newblock Preprint arXiv:1408.3092.

\bibitem{tmvtnorm}
S. Wilhelm,
\newblock Package ``tmvtnorm'',
\newblock {\em http://cran.r-project.org/web/packages/tmvtnorm/}

\bibitem{zhou2010nonparametric}
M.~Zhou, C.~Wang, M.~Chen, J.~Paisley, D.~Dunson, and L.~Carin.
\newblock Nonparametric bayesian matrix completion.
\newblock {\em Proc. IEEE SAM}, 2010.

\end{thebibliography}

\end{document}